\documentclass[12pt]{article}
\usepackage{amssymb}
\usepackage{latexsym}
\newtheorem{Th}{Theorem}
\newtheorem{Prop}[Th]{Proposition}

\newtheorem{Rem}[Th]{\it Remark}
\newtheorem{Cor}[Th]{Corollary}

\newcommand{\be}{\begin{eqnarray*}}
\newcommand{\ee}{\end{eqnarray*}}

\newcommand{\pkE}{{\cal P}(^k\!E)}

\newcommand{\pol}{{\cal P}(}
\newcommand{\WCo}{{\cal WC}o}

\newcommand{\Proof}{\noindent {\it Proof. }}
\newcommand{\fin}{\hspace*{\fill} $\Box$}
\newcommand{\finesp}{\hspace*{\fill} $\Box$\vspace{.5\baselineskip}}
\newcommand{\R}{{\if mm {\rm I}\mkern -3mu{\rm R}\else \leavevmode
        \hbox{I}\kern -.17em \hbox{R} \fi}}
\newcommand{\K}{{\if mm {\rm I}\mkern -3mu{\rm K}\else \leavevmode
        \hbox{I}\kern -.17em \hbox{K} \fi}}

\newcommand{\C}{{\if mm {{\rm C}\mkern -15mu{\phantom{\rm t}\vrule}}
     \mkern +10mu \else \leavemode \hbox{I}\kern -.17em \hbox{C} \fi}}
\newcommand{\Q}{{\if mm {{\rm Q}\mkern -16mu{\phantom{\rm t}\vrule}}
     \mkern +10mu \else \leavemode \hbox{I}\kern -.17em \hbox{Q} \fi}}
\newcommand{\N}{\mathbb{N}}
\newcommand{\ra}{\rightarrow}
\newcommand{\lra}{\longrightarrow}
\newcommand{\Ra}{\Rightarrow}       
\newcommand{\espv}{\vspace{.5\baselineskip}}
\newcommand{\espvv}{\vspace{\baselineskip}}
\newcommand{\noin}{\noindent}

\newcommand{\eps}{\epsilon}
\def\card{\mathop{\rm card}}
\def\supp{\mathop{\rm supp}}

\begin{document}
\title{Polynomials on Schreier's space}

\author{Manuel Gonz\'alez\thanks{Supported in part by DGICYT Grant PB
      97--0349 (Spain)}        \and
      Joaqu\'\i n  M. Guti\'errez\thanks{Supported in part by DGICYT
      Grant PB 96--0607 (Spain)}}

\date{\hspace*{\fill}}
\maketitle

\vspace{-1\baselineskip}

\begin{abstract}
We introduce a weakened version of the Dunford-Pettis
property, and give examples of Banach spaces with this property.
In particular, we
show that every closed subspace of
Schreier's space $S$ enjoys it.
As an application, we characterize the weak polynomial convergence
of sequences, show that every closed subspace of $S$ has the
polynomial Dunford-Pettis property of Bistr\"om et al.\
and give other polynomial properties of $S$.
\end{abstract}


\noin
\hrule\espvv

A subset $A=\{ n_1<\cdots <n_k\}$ of the natural numbers $\N$ is said
to be {\it admissible\/} if $k\leq n_1$.
Schreier's space $S$ \cite{Sc,CS} is the completion of the space
$c_{00}$ of all
scalar sequences of finite support with respect to the norm:
$$
\|x\|_S:=\sup\left\{\sum_{j\in A}|x_j|: A\subset\N \mbox{ is admissible}
\right\}\, ,\hspace{1em}\mbox{ for } x=\left( x_j\right)_{j=1}^\infty
\, .
$$
Some basic properties of $S$ may be seen in \cite{CG}. Schreier's
space has been used to provide counterexamples in Banach space
theory \cite{Be,CG,CGIs,Os,PS}.

In this paper, we introduce a weakened version of the
Dunford-Pettis property, and give examples of Banach spaces with
this property. In particular, we show that every closed subspace
of $S$ enjoys it. It is well known that a reflexive Banach space
with the Dunford-Pettis property must be finite dimensional. The
same is true for a Banach space with the Banach-Saks property and
the weak Dunford-Pettis property. As an application, we
investigate polynomial properties of $S$, characterizing the
sequences which converge in the weak polynomial topology, that we
shall call the $\cal P$-topology. As far as we know, this is the
first time that $\cal P$-convergent sequences are characterized
for a space where $\cal P$-convergence does not coincide with
either norm or weak convergence of sequences. From this we obtain
that every closed subspace of $S$ has the polynomial
Dunford-Pettis property \cite{BJL}.

We also show that the relatively compact sets for the $\cal P$-topology
coincide with the Banach-Saks sets, that the absolutely convex
closed hull of a Banach-Saks set in $S$ is a Banach-Saks set,
and that the tensor product of two Banach-Saks sets is a Banach-Saks
set in the projective tensor product $S\otimes_\pi S$.
It is unknown if the Banach-Saks sets in an arbitrary Banach space
are stable under convex hulls.
An example of a Banach space so that the relatively
$\cal P$-compact sets are not stable under convex hulls was given in
\cite{CGG}.
Moreover, given two $\cal P$-null
(i.e., $\cal P$-convergent to zero) sequences $(x_n), (y_n)\subset S$,
we prove that $\{ x_n\otimes y_n\}$ is a Banach-Saks set in
$S\otimes_\pi S$.
The polynomial Dunford-Pettis property of $S$ implies that the sequence
$(x_n\otimes y_n)$ is $\cal P$-null
in $S\otimes_\pi S$, and that $(x_n+y_n)$ is $\cal P$-null in $S$.
These properties have interesting consequences
in infinite dimensional holomorphy, as shown in \cite[Remark~4.7]{Gant}.

We shall use the facts that the unit vector basis of $S$ is
unconditional, and that every
closed subspace of $S$ contains an isomorphic copy of $c_0$
(so, $S$ contains no copy of $\ell_1$).\espv

Throughout the paper, $E$ will denote a Banach space, and $E^*$
its dual. The space of all scalar valued $k$-homogeneous
(continuous) polynomials on $E$ is represented by $\pkE$. General
references for polynomials on Banach spaces are \cite{Din,Mu}.
Given a subset $A\subset\N$, $\card A$ stands for the cardinality
of $A$.

A sequence $(x_n)\subset E$ is {\it $\cal P$-convergent\/} to $x$ if
$P(x_n)\ra P(x)$ for every $P\in \pkE$ and all $k\in\N$.
A set $A\subset E$ is {\it relatively $\cal P$-compact\/} if every
sequence in $A$ has a $\cal P$-convergent subsequence.

A subset $A\subset E$ is a {\it Banach-Saks set\/} if every
sequence in $A$ has a subsequence whose arithmetic means converge
in norm. A sequence $(x_n)\subset E$ {\it converges uniformly
weakly\/} to $x$ in $E$ \cite[Definition~2.1]{Me} if, for each
$\eps>0$, there exists $N(\eps)\in\N$ such that $\card\{ n\in\N:
|\phi (x_n-x)|\geq \eps\}\leq N(\eps)$, for every $\phi\in E^*$
with $\|\phi\|\leq 1$. A subset $A\subset E$ is a Banach-Saks set
if and only if every sequence in $A$ has a subsequence which is
uniformly weakly convergent in $E$ \cite[Theorem~2.9]{Me}.

Recall that
a Banach space $E$ has the {\it Dunford-Pettis property\/}
(DPP, for short) if, for all weakly null
sequences $(x_n)\subset E$ and
$(\phi_n)\subset E^*$, we have $\phi_n(x_n)\ra 0$.
We say that $E$ has the {\it polynomial Dunford-Pettis property\/}
if, for every ${\cal P}$-null sequence $(x_n)\subset E$ and every
weakly null sequence $(\phi_n)\subset E^*$, we have $\phi_n(x_n)\ra 0$.
The DPP implies the polynomial DPP.
$E$ is said to be a
{\it $\Lambda$-space\/} if ${\cal P}$-null sequences
and norm null sequences coincide in $E$.
Spaces with the Schur property are trivially $\Lambda$-spaces.
All superreflexive spaces are $\Lambda$-spaces \cite{JP}.
It is proved in \cite[Corollary~3.6]{FJ} that every Banach space
with nontrivial type is a $\Lambda$-space.

\section{The weak Dunford-Pettis property}

We say that a Banach space $E$ has the {\it weak Dunford-Pettis
property\/} (wDPP, for short) if, given a uniformly weakly null
sequence $(x_n)\subset E$ and a weakly null sequence
$(\phi_n)\subset E^*$, we have $\lim\phi_n(x_n)=0$.

The space $\ell_2$ fails the wDPP since its unit vector basis is
uniformly weakly null. Clearly, if $E$ has the DPP, then $E$ has
the wDPP.

Denote by $T$ the dual of the original Tsirelson space $T^*$
\cite{CS}. Then the uniformly weakly convergent sequences in $T$
are norm convergent. Indeed, suppose $(x_n)$ is uniformly weakly
convergent to $x\in T$ and $\|x_n-x\|\geq\delta>0$. Passing to a
subsequence, we may assume that the sequence $(x_n-x)$ is basic
and equivalent to a subsequence of the unit vector basis $(t_n)$
of $T$
\cite[Proposition~II.7]{CS}. If $A\subset\N$ is admissible, by the
definition of the norm of $T$, we have:
$$
\left\|\sum_{i\in A}t_i\right\|\geq\frac{1}{2}\card A
$$
and so, $(t_n)$ has no uniformly weakly null subsequence,
which yields a contradiction.

Therefore, $T$ enjoys the wDPP, but $T^*$ does not since the unit
vector basis of $T^*$ is a Banach-Saks set. We conclude that the
wDPP of a Banach space neither implies nor is implied by the wDPP
of its dual.

The following simple remark will be useful:

\begin{Prop}
A Banach space $E$ has the {\rm wDPP} if and only if whenever
$(x_n)\subset E$ is uniformly weakly null and $(\phi_n)\subset
E^*$ is weak Cauchy, we have $\lim\phi_n(x_n)=0$.
\end{Prop}

\Proof
For the nontrivial part, if $\phi_n(x_n)\geq\delta>0$, we can find
$k_1<\cdots<k_n<\cdots$ such that $\left|\phi_n\left(
x_{k_n}\right)\right|<\delta/2$. Then,
$$
\delta\leq\phi_{k_n}\left( x_{k_n}\right)\leq\left|\left(\phi_{k_n}
-\phi_n\right)\left( x_{k_n}\right)\right| +\left|\phi_n\left(
x_{k_n}\right)\right|
$$
and the right hand side is less than $\delta$ for $n$ large enough,
since the sequence $\left(\phi_{k_n}-\phi_n\right)$ is weakly
null.\finesp

Denoting by $\WCo (E,F)$ the space of all weakly compact (linear)
operators from $E$ into the Banach space $F$, and by ${\cal C}_w(E,F)$
the space of all operators taking uniformly weakly null sequences
in $E$ into norm null sequences in $F$, we have:

\begin{Prop}
The Banach space $E$ satisfies the {\rm wDPP} if and only if, for all
Banach spaces $F$, we have $\WCo (E,F)\subseteq {\cal C}_w(E,F)$.
\end{Prop}

\Proof
Suppose $E$ has the wDPP and $(x_n)\subset E$ is uniformly weakly null.
Take $L\in\WCo (E,F)$ with adjoint $L^*$.
Choose $(\phi_n)$ in the unit ball of $F^*$ such that
$\phi_n(Lx_n)=\|Lx_n\|$. There is a subsequence
$\left(\phi_{n_k}\right)$ such that $\left(L^*\phi_{n_k}\right)$ is
weakly convergent. Hence, $\phi_n(Lx_n)=(L^*\phi_n)x_n\ra 0$.
Conversely, if $E$ fails the wDPP, we can find $(x_n)$ uniformly weakly
null in $E$ and $(\phi_n)$ weakly null in $E^*$ such that
$\phi_n(x_n)\geq\delta>0$. We define an operator $L:E\ra c_0$ by
$Lx:=(\phi_n(x))$. Then, $L$ is weakly compact but $\|Lx_n\|\geq|\phi_n
(x_n)|\geq\delta>0$ for all $n$.\finesp

The following easy fact characterizes the reflexive Banach spaces with
the wDPP:

\begin{Prop}
Let $E$ be a reflexive Banach space. Then $E$ has the {\rm wDPP}
if and only if every uniformly weakly null sequence in $E$ is norm null.
\end{Prop}

\Proof
Suppose there is a uniformly weakly null sequence $(x_n)\subset E$
with $\|x_n\|=1$. We can assume that $(x_n)$ is basic and the
sequence of coefficient functionals $(\phi_n)$ is weakly null in
$E^*$. Since $\phi_n(x_n)=1$, we conclude that $E$ does not have
the wDPP. The converse is clear.\finesp

Recall that a Banach space $E$ has the {\it Banach-Saks property\/} if
every bounded subset in $E$ is a Banach-Saks set. We then have:

\begin{Cor}
If $E$ has the Banach-Saks property and the {\rm wDPP}, then $E$ is
finite dimensional.
\end{Cor}

A space $E$ has the {\it weak Banach-Saks property\/} if every
weakly null sequence in $E$ contains a subsequence whose
arithmetic means converge. Equivalently \cite[Theorem~2.10]{Me},
every weakly null sequence has a subsequence which converges to
zero uniformly weakly in $E$. The space $L^1[0,1]$ has the weak
Banach-Saks property. The following result is clear:

\begin{Prop}
Assume $E$ has the weak Banach-Saks property. Then $E$ has the
{\rm DPP} if and only if $E$ has the {\rm wDPP}.
\end{Prop}

We say that $E$ has the {\it hereditary weak Dunford-Pettis
property\/} if every closed subspace of $E$ has the wDPP.

\begin{Prop}
\label{nuwns}
A Banach space $E$ has the hereditary {\rm wDPP} if and only if
every normalized uniformly weakly null sequence in $E$ contains a
subsequence equivalent to the $c_0$-basis.
\end{Prop}

\Proof
Suppose that the uniformly weakly null sequence $(x_n)\subset E$,
$\|x_n\|=1$, has no subsequence equivalent to the $c_0$-basis. We
can assume that $(x_n)$ is basic. Let $(\phi_n)\subset [x_n]^*$ be
the sequence of coefficient functionals, where $[x_n]$ denotes the
closed linear span of the set $\{ x_n\}$ in $E$. After taking a
subsequence, we can assume that either $(\phi_n)$ is equivalent to
the $\ell_1$-basis or $(\phi_n)$ is weak Cauchy \cite[Ch.~XI]{Di}.
In the first case, we define an operator $L:[x_n]\ra c_0$ by
$L(x):=(\phi_n(x))$. Clearly, $L$ is injective and has dense
range. The adjoint $L^*:\ell_1\ra [x_n]^*$ takes the unit vector
basis of $\ell_1$ into the sequence $(\phi_n)$ and has therefore
closed range. Hence, $L$ is a surjective isomorphism, which
contradicts our assumption. So, $(\phi_n)$ must be weak Cauchy.
Since $\phi_n(x_n)=1$, the subspace $[x_n]$ fails to have the
wDPP.

For the converse, it is enough to show that $E$ has the wDPP.
Suppose it does not. Then we can find a uniformly weakly null
sequence $(x_n)\subset E$ and a weakly null sequence
$(\phi_n)\subset E^*$ such that $\phi_n(x_n)\geq 1$ for all $n$.
Passing to a subsequence, we can assume that $(x_n)$ is equivalent
to the $c_0$-basis. Since the dual of $c_0$ has the Schur
property, the restriction of $(\phi_n)$ to the subspace $[x_k]$ is
norm null, and we get a contradiction.\fin

\begin{Rem}
{\rm This simple proof also shows that a Banach space $E$ has the
hereditary {\rm DPP} if and only if every normalized weakly null
sequence in $E$ has a subsequence equivalent to the $c_0$-basis
\cite[Proposition~2]{Ce}. From this we get that every infinite
dimensional Banach space without a copy of either $c_0$ or
$\ell_1$ contains a subspace without the DPP \cite[p.~254]{Di}.
The original proofs of these two results were based on a
characterization of $c_0$'s unit vector basis that Elton \cite{El}
obtained by using Ramsey's theorem. }
\end{Rem}

Our aim now is to show that Schreier's space enjoys the hereditary
wDPP.

\begin{Prop}
\label{infnorm}
If $(x_n)$ is a uniformly weakly null sequence in $S$, then
$\|x_n\|_\infty\ra 0$.
\end{Prop}

\Proof
Let $x_n=\left( x_n^i\right)_{i=1}^\infty$.
Since a set of $\pm 1$'s on an admissible set
is a norm-one functional on $S$,
given $\eps>0$, there is $N(\eps)\in\N$ such that
$$
\card\left\{ n\in\N :\sum_{i\in A}\left| x_n^i\right|\geq \eps\right\}
\leq N(\eps)
$$
for each admissible $A$.
Suppose our statement fails; then we can find $\delta>0$ and two
increasing sequences of indices $(n_k), (l_k)$ such that
$$
\left| x_{n_k}^{l_k}\right|\geq \delta\hspace{1em}\mbox{ for all $k$}.
$$
The set $A_m:=\{ l_{m+1},\ldots,l_{2m}\}$ is admissible for each
$m\in\N$, and
$$
\card\left\{ n\in\N :\sum_{i\in A_m}\left|
x_n^i\right|\geq\delta\right\}\geq m\, ,
$$
a contradiction which finishes the proof.\finesp

The converse is not true. Indeed, take
$x_n:=(e_1+\cdots+e_n)/n$.
The set $A_k:=\{ 2^{k-1},\ldots,2^k-1\}$ is admissible for each
$k\in\N$.
Denoting by $\left( e^*_i\right)$ the unit vector basis of $S^*$,
the functional
$$
\phi_k:=\sum_{i=2^{k-1}}^{2^k-1} e^*_i\in S^*
$$
has norm one.
Choosing $n$ so that $2^{k-2}+2^{k-1}\leq n\leq 2^k-1$, we have
$$
\phi_k(x_n)\geq\frac{2^{k-2}}{n}>\frac{2^{k-2}}{2^k}
=\frac{1}{4}\; .
$$
Therefore, $\|x_n\|_\infty\ra 0$, but $(x_n)$ does not converge to
zero uniformly weakly. The proof of the following result is
essentially contained in
\cite{CGIs}. We give it for completeness.

\begin{Prop}
\label{c0basis}
Let $(x_n)$ be a normalized sequence in $S$ such that
$\|x_n\|_\infty\ra 0$. Then $(x_n)$ contains a subsequence
equivalent to the $c_0$-basis.
\end{Prop}

\Proof
Let us denote by $\supp (x)$ the support of $x$. Passing to a
subsequence and perturbing it with a null sequence, we can assume
that $\max\supp (x_n)<\min\supp (x_{n+1})$, and
\begin{eqnarray}
\label{cota}
\|x_n\|_\infty\leq\frac{1}{2^n\max\supp (x_{n-1})}\; .
\end{eqnarray}
Given $x_{n_1},\ldots,x_{n_m}$ and an admissible set $A$, we take
$k_0$ to be the minimum value of $k$ such that
$A\cap\supp\left( x_{n_k}
\right)\neq\emptyset$.
In particular, this implies that
$\card A\leq \max\supp\left( x_{n_{k_0}}\right)$.
Denoting $x_n(i):=x_n^i$, we have:
\be
\sum_{i\in A}\left|\left(\sum_{k=1}^m x_{n_k}\right)(i)\right|&=&
\sum_{i\in A}\left|\left(\sum_{k=k_0}^m x_{n_k}\right)(i)\right|\\
&=&\sum_{k=k_0}^m\sum_{i\in A\cap\supp\left( x_{n_k}\right)}
   \left| x_{n_k}(i)\right|\\
&\leq&\left\| x_{n_{k_0}}\right\| + \sum_{k=k_0+1}^m
      \left\| x_{n_k}\right\|_\infty\cdot\card A\\
&\leq&\left\| x_{n_{k_0}}\right\| + \sum_{k=k_0+1}^m2^{-n_k}\\
&\leq&2\, ,
\ee
where we have used (\ref{cota}).
Thus we have proved that
$$
\left\|\sum_{k=1}^mx_{n_k}\right\|\leq 2
$$
and hence the series $\sum x_n$ is weakly unconditionally Cauchy.
Therefore, $(x_n)$ has a subsequence equivalent to the $c_0$-basis
\cite[Corollary~V.7]{Di}.\finesp

Combining the last two results with Proposition~\ref{nuwns} yields:

\begin{Th}
Schreier's space $S$ has the hereditary {\rm wDPP}.
\end{Th}

We now show that the dual $S^*$ of Schreier's space fails the wDPP.
The next result follows the lines of \cite[Example~13,(H)]{Le}.

\begin{Prop}
\label{l1basis}
Let $(\phi_n)$ be a normalized block basis of the unit basis of
$S^*$ such that $\|\phi_n\|_\infty\ra 0$. Then $(\phi_n)$ contains
a subsequence equivalent to the $\ell_1$-basis.
\end{Prop}

\Proof
Let $(x_n)$ be a sequence in $S$ such that $\|x_n\| < 2$,
$\supp (x_n) = \supp(\phi_n)$ and $\phi_n(x_n) = 1$ for every $n$.

First we select $n_1$ such that $\min\supp\left(\phi_{n_1}\right)
> 2^2$ and $\left\|\phi_{n_1} \right\|_{\infty} < 2^{-4}$. Since
$\left\|x_{n_1}\right\| < 2$, the set
$$
A_1 = \left\{ i \in \N : \left|x_{n_1}(i)\right| \geq 2^{-1} \right\}
$$
has fewer than $2^2$ elements. We define $y_{n_1}(i) = 0$ if $i\in
A_1$ and $y_{n_1}(i) = x_{n_1}(i)$ otherwise, and obtain $y_{n_1}
\in S$ such that $\left\|y_{n_1}\right\| < 2$,
$\left\|y_{n_1}\right\|_{\infty} < 2^{-1}$ and
$$
\left|\phi_{n_1}\left(y_{n_1}\right)\right| \geq \phi_{n_1}(x_{n_1}) -
\left|\phi_{n_1}(y_{n_1} - x_{n_1})\right| > 1 - 2 (2^2) 2^{-4}
= 2^{-1}.
$$

Next we select $n_2 > n_1$ such that $\min\supp\left(\phi_{n_2}\right)
> 2^3$ and $\left\|\phi_{n_2} \right\|_{\infty} < 2^{-5}$.
Since $\left\|x_{n_2}\right\| < 2$, the set
$$
A_2 = \left\{ i \in \N : \left|x_{n_2}(i)\right| \geq 2^{-2} \right\}
$$
has fewer than $2^3$ elements. We define $y_{n_2}(i) = 0$ if $i\in
A_2$ and $y_{n_2}(i) = x_{n_1}(i)$ otherwise, and obtain $y_{n_2}
\in S$ such that $\left\|y_{n_2}\right\| < 2$,
$\left\|y_{n_2}\right\|_{\infty} < 2^{-2}$ and
$$
\left|\phi_{n_2}\left(y_{n_2}\right)\right| \geq \phi_{n_2}(x_{n_2}) -
\left|\phi_{n_2}(y_{n_2} - x_{n_2})\right| > 1 - 2(2^3) 2^{-5}
= 2^{-1}.
$$

In this way we get a subsequence $\left(\phi_{n_j}\right)$ and
a sequence $\left( y_{n_j}\right)\subset S$ such that
$\left|\phi_{n_j}\left( y_{n_j}\right)\right| > 2^{-1}$,
$\left\|y_{n_j}\right\| < 2$ and $\left\|y_{n_j}\right\|_{\infty}
< 2^{-j}$.
Passing to a subsequence we can assume by Proposition~\ref{c0basis}
that $\left( y_{n_j}\right)$ is equivalent to the $c_0$-basis,
from which it easily follows that $\left(\phi_{n_j}\right)$
is equivalent to the $\ell_1$-basis.\fin

\begin{Prop}
\label{weakBS}
The dual $S^*$ of Schreier's space $S$ has the weak Banach-Saks
property.
\end{Prop}

\Proof
Let $(\phi_n)$ be a normalized weakly null sequence in $S^*$.
Passing to a subsequence we can assume that $(\phi_n)$ is
equivalent to a block basis of the unit basis. We have that
$\left((\phi_1+\cdots+\phi_n)/n)\right)$ is a weakly null sequence
and $\left\|(\phi_1+\cdots+\phi_n)/n\right\|_{\infty}\ra 0$. If
$\left\|(\phi_1+\cdots+\phi_n)/n\right\|$ does not converge to
$0$, passing to a subsequence, it follows from
Proposition~\ref{l1basis} that
$\left((\phi_1+\cdots+\phi_n)/n)\right)$ contains a subsequence
equivalent to the $\ell_1$-basis, a contradiction.\fin

\begin{Cor}
The dual $S^*$ of Schreier's space does not have the {\rm wDPP}.
\end{Cor}

\section{Applications to polynomials}

In this Section, we describe the $\cal P$-convergence of sequences
in $S$, obtaining thereby some polynomial properties of this
space, and characterize the Banach-Saks sets in it.

We shall use the fact that $S$ may be algebraically embedded in
$\ell_2$, and that the natural inclusion $j:S\ra\ell_2$ is
continuous. To see this, take $x:=(x_i)\in S$, $\|x\|_S=1$, and
call $y:=(y_i)$ the sequence $(|x_i|)$, reordered in a
nonincreasing way. Then $\|y\|_2=\|x\|_2$, and $\|y\|_S\leq 1$.
This implies $y_{2k-1}\leq k^{-1}$ for each $k$. Therefore,
$$
\|y\|_2^2=\sum_{i=1}^\infty y_i^2\leq
1+1+\frac{1}{2^2}+\frac{1}{2^2}+\frac{1}{3^2}+\frac{1}{3^2}+\cdots
=\frac{\pi^2}{3}\; ,
$$
from which $\|j\|\leq\pi/\sqrt{3}$.

As a consequence, $P(x):=\|x\|_2^2$ defines a
$2$-homogeneous polynomial on $S$.

\begin{Prop}
\label{pnull}
Let $(x_n)$ be a sequence in $S$. The following assertions are
equivalent:

{\rm (a)} $(x_n)$ is $\cal P$-null;

{\rm (b)} $(x_n)$ is bounded in $S$ and $\|x_n\|_2\ra 0$;

{\rm (c)} $(x_n)$ is bounded in $S$ and $\|x_n\|_\infty\ra 0$.
\end{Prop}

\Proof
(a) $\Ra$ (b) since $P(x):=\|x\|^2_2$ is a polynomial on $S$.

(b) $\Ra$ (c) is clear.

(c) $\Ra$ (a). It is enough to show that $(x_n)$ has a $\cal
P$-null subsequence. If $\inf\|x_n\|>0$, then there is a
subsequence of $(x_n)$ equivalent to the $c_0$-basis
(Proposition~\ref{c0basis}) and so $\cal P$-null, since the
$c_0$-basis is $\cal P$-null. If $\inf\|x_n\|=0$, then there is a
norm null subsequence, which is $\cal P$-null a fortiori.\finesp

A Banach space has the {\it hereditary polynomial\/} DPP if
every closed subspace has the polynomial DPP.

\begin{Th}
The space $S$ has the hereditary polynomial {\rm DPP}.
\end{Th}

\Proof
By Propositions~\ref{pnull} and~\ref{c0basis}, every normalized
${\cal P}$-null sequence in $S$ contains a subsequence equivalent
to the $c_0$-basis. Obvious modifications in the ``if" part
of the proof of Proposition~\ref{nuwns} yield the result.\finesp

It is shown in \cite{BJL} that, given two ${\cal P}$-null sequences
$(x_n)$, $(y_n)$ in a space with the polynomial DPP, the sequence
$(x_n+y_n)$ is ${\cal P}$-null.
A Banach space where this is not true was recently found by Castillo
et al.\ \cite{CGG}.

\begin{Prop}
Let $A$ be a subset of $S$. The following assertions are equivalent:

{\rm (a)} $A$ is a Banach-Saks set;

{\rm (b)} $A$ is relatively $\cal P$-compact;

{\rm (c)} $A$ is relatively weakly compact in $S$, and relatively
compact as a subset of $\ell_\infty$.
\end{Prop}

\Proof
(a) $\Ra$ (b). Let $A$ be a Banach-Saks set. Given a sequence
$(x_n)\subset A$, passing to a subsequence, we may assume that $(x_n)$
converges to some $x$ uniformly weakly in $S$. Then, $(x_n-x)$
has a subsequence which is either
norm null or equivalent to the $c_0$-basis.
In both cases, $(x_n)$ is $\cal P$-convergent to $x$.

(b) $\Ra$ (c). If $A$ is relatively $\cal P$-compact, it is
relatively weakly compact. Moreover, given a sequence
$(x_n)\subset A$, we can assume that $(x_n-x)$ is $\cal P$-null
for some $x$. By Proposition~\ref{pnull}, $\|x_n-x\|_\infty\ra 0$
and so $A$ is relatively compact as a subset of $\ell_\infty$.

(c) $\Ra$ (a). Choose a sequence $(x_n)\subset A$. We may assume that
$(x_n)$ is weakly convergent to some $x$ and $\|x_n-x\|_\infty\ra 0$.
Passing to a subsequence, we have either $\|x_n-x\|\ra 0$ or,
by Proposition~\ref{c0basis}, $(x_n-x)$ is equivalent to
the $c_0$-basis and is therefore uniformly weakly null.\fin

\begin{Cor}
If $A$ is a Banach-Saks set in $S$, then the absolutely convex closed
hull of $A$ is a Banach-Saks set.
\end{Cor}

The following two properties were introduced in \cite{ACL} and
studied by various authors (see, e.g., \cite{BJL,CK}).

(a) A Banach space $E$ has property (P) if, given two bounded
sequences $(u_n),(v_n)$ in $E$ such that $P(u_n)-P(v_n)\ra 0$ for
every $P\in\pkE$ and all $k$, it follows that the sequence
$(u_n-v_n)$ is $\cal P$-null. Every superreflexive space and every
space with the DPP have property (P). A Banach space failing to
have property (P) has been found by Castillo et al.\ \cite{CGG}.

(b) A Banach space $E$ has property (RP) if, given two bounded
sequences $(u_n),(v_n)$ in $E$ such that the sequence $(u_n-v_n)$
is $\cal P$-null, it follows that $P(u_n)-P(v_n)\ra 0$ for every
$P\in\pkE$ and all $k$. Every $\Lambda$-space and every predual of
a Banach space with the Schur property have property (RP). The
spaces $L_1[0,1]$, $C[0,1]$ and $L_\infty [0,1]$ fail to have
property (RP) \cite{ACL}.

We now show that $S$ has property (P) and fails property (RP).

\begin{Prop}
The space $S$ fails property {\rm (RP).}
\end{Prop}

\Proof
Consider the vectors
$$
v_n:=e_n\; ;\hspace{1cm}u_n:=e_n+2^{1-n}\left( e_{2^{n-1}}+\cdots
+e_{2^n-1}\right)\, .
$$
Then $\|u_n-v_n\|_\infty\ra 0$ and so $(u_n-v_n)$ is $\cal P$-null in
$S$.
Define
$$
P(x):= \sum_{n=1}^\infty x_n^2\left(\sum_{k=2^{n-1}}^{2^n-1}x_k
\right)\; ,\hspace{1em}\mbox{ for $x=(x_n)\in S$}.
$$
Since
$$
|P(x)|\leq\|x\|_S\cdot\|x\|_2^2\leq \frac{\pi^2}{3}\cdot\|x\|_S^3\, ,
$$
we get that $P\in\pol ^3\! S)$.
We have $P(v_n)=0$ and $P(u_n)=1$ for all $n>1$.\finesp

In the above proof, we need a polynomial of degree greater than or
equal to three. Indeed, if $P\in\pol ^2\! S)$ and
$(u_n),(v_n)\subset S$ are bounded with $(u_n-v_n)$ $\cal P$-null,
denoting $w_n:=u_n-v_n$, we have
$$
P(u_n)-P(v_n)=P(w_n+v_n)-P(v_n)=2\hat{P}(w_n,v_n)+P(w_n)\, ,
$$
where $\hat{P}$ is the symmetric bilinear form associated to $P$.
Let $\bar{P}:S\ra S^*$ be the operator defined by
$\bar{P}(x)(y):=\hat{P}(x,y)$. Since $S$ has an unconditional
basis and contains no copy of $\ell_1$, the space $S^*$ has an
unconditional basis and is weakly sequentially complete.
Therefore, every operator from $S$ into $S^*$ is weakly compact.
Passing to a subsequence, we can assume that $(w_n)$ is uniformly
weakly null. Since $S$ has the wDPP, $\|\bar{P}(w_n)\|\ra 0$.
Hence, $\hat{P}(w_n,v_n)=\bar{P}(w_n)(v_n)\ra 0$. Clearly,
$P(w_n)\ra 0$ and so $P(u_n)-P(v_n)\ra 0$.

\begin{Prop}
The space $S$ enjoys property {\rm (P)}.
\end{Prop}

\Proof
Let $(u_n),(v_n)\subset S$ be bounded sequences such that
$(u_n-v_n)$ is not $\cal P$-null. We wish to find $Q\in\pol
^k\!S)$ for some $k$ so that $(Q(u_n)-Q(v_n))$ does not tend to
zero. By $u_n^i$ and $v_n^i$ we shall denote the $i\,$th
coordinate of $u_n$ and $v_n$, respectively.

If $(u_n-v_n)$ is not weakly null, then
$\phi(u_n)-\phi(v_n)\not\ra 0$ for some $\phi\in S^*$. It is
enough to take $Q:=\phi$.

If $(u_n-v_n)$ is weakly null, passing to a subsequence and
perturbing it by a norm null sequence, we can assume that
$(u_n-v_n)$ is a block basis:
$$
u_n-v_n=\sum_{i=k_n}^{l_n}a_ie_i\, .
$$
Take $p_n$ with $k_n\leq p_n\leq l_n$ and $\left| a_{p_n}\right|
=\|u_n-v_n\|_\infty$.
We know that $\|u_n-v_n\|_\infty$ does not go to zero.
Passing to a subsequence, we may assume:
$$
v_n^{p_n}\lra v\, ;
\hspace{1em}u_n^{p_n}\lra u\, ;\hspace{1em}u\neq v\, .
$$
Let $P\in\pol ^2\! S)$ be given by
$$
P(x):=\sum_{i=1}^\infty \left( x_{p_i}\right)^2\, .
$$
If $P(u_n)-P(v_n)\not\ra 0$, we are done.
If
$$
0=\lim\left[ P(u_n)-P(v_n)\right] = \lim\left[\left( u_n^{p_n}\right)
^2-\left( v_n^{p_n}\right)^2\right] =u^2-v^2=(u-v)(u+v)\, ,
$$
we have $u=-v=\alpha$, for some $\alpha\neq 0$.
Defining $Q\in\pol ^3\!S)$ by
$$
Q(x):=\sum_{i=1}^\infty \left(x_{p_i}\right)^3\, ,
$$
we have $Q(u_n)-Q(v_n)\ra 2\alpha^3\neq 0$,
and the proof is finished.\fin

\begin{Prop}
Let $(x_n),(y_n)\subset S$ be $\cal P$-null sequences. Then:

{\rm (a)}
the set $\{ x_n\otimes y_n\}$ is a Banach-Saks set in
$S\otimes_\pi S$;

{\rm (b)}
the sequence $( x_n\otimes y_n)$ is $\cal P$-null in $S\otimes_\pi S$.
\end{Prop}

\Proof
(a) Since $(x_n)$ and $(y_n)$ have subsequences equivalent to the
$c_0$-basis, it is enough to show that $(e_n\otimes e_n)$ is
uniformly weakly null in $c_0\otimes_\pi c_0$. Take $L\in
(c_0\otimes_\pi c_0)^*$, which may be viewed as an operator from
$c_0$ into $\ell_1$. Since the series $\sum e_n$ is weakly
unconditionally Cauchy, using \cite[Theorem~2]{GGUC}, we can find
$C>0$ such that $\sum|\langle Le_n,e_n\rangle|\leq C$ whenever
$\|L\|\leq 1$. Therefore, given $\eps>0$, choosing $N\in\N$ with
$N\geq C/\eps$, we have
$$
\card\left\{ n\in\N :|\langle Le_n,e_n\rangle|\geq\eps\right\}\leq N
$$
if $\|L\|\leq 1$, and the result is proved.

(b) Since $S$ has the polynomial DPP, part (b) follows from
\cite[Theorem~2.1]{BJL}.\finesp

As a consequence, if $A,B\subset S$ are Banach-Saks sets, then
$A\otimes B$ is a Banach-Saks set in $S\otimes_\pi S$.\espv

{\small
\begin{tabular}{ll}
Manuel Gonz\'alez               &Joaqu\'\i n M. Guti\'errez\\
Departamento de Matem\'aticas   &Departamento de Matem\'aticas\\
Facultad de Ciencias            &ETS de Ingenieros Industriales\\
Universidad de Cantabria        &Universidad Polit\'ecnica de Madrid\\
39071 Santander (Spain)         &C. Jos\'e Guti\'errez Abascal 2\\
e-mail gonzalem@ccaix3.unican.es\phantom{unican}&28006 Madrid (Spain)\\
                                &e-mail gutierrezj@member.ams.org
\end{tabular}}

\vspace{\fill}\hspace*{\fill}{\small file poss.tex}

\end{document}